\pdfoutput=1                            
\def\lcversion{}					    

\documentclass[nonote,noversion,squeezed]{cdcarticle}

\def\MODE{2}

\usepackage{amsmath}
\usepackage{mathrsfs}
\usepackage{multirow}
\usepackage{textcomp}
\usepackage{graphicx} 
\usepackage{tabularx}
\usepackage{url}
\usepackage{array,tabularx}

\usepackage[utf8]{inputenc}				
\usepackage[english]{babel}				

\usepackage{math}						
\usepackage{enumitem}
\usepackage{cite}
\usepackage[font=small,labelfont=bf]{caption}

\usepackage[hidelinks]{hyperref}
\usepackage{subcaption}

\usepackage{mdframed}
\newmdtheoremenv[innertopmargin=0,innerbottommargin=5,%
innerleftmargin=5,innerrightmargin=5]{fthm}[thm]{Theorem}

\usepackage{tikz,pgfplots}
\usetikzlibrary{external,calc,positioning}
\usepgfplotslibrary{fillbetween}
\usepgfplotslibrary{external,colorbrewer}
\pgfplotsset{compat=1.14}
\graphicspath{{./graphics/}}

\let\oldbibliography\thebibliography
\renewcommand{\thebibliography}[1]{\oldbibliography{#1}
\setlength{\itemsep}{1pt}} 

\usepackage{booktabs}
\newcommand{\bbmat}[1]{{\renewcommand{\arraystretch}{1.0}\addtolength{\arraycolsep}{-1mm}\bmat{#1}}}

\newcommand{\matbox}[1]{\begin{minipage}[m]{1cm}\centering$\bbmat{#1}$\\[-0.5\baselineskip]\end{minipage}}

\begin{document}

\title{A Revised Mehrotra Predictor-Corrector algorithm for Model Predictive Control}

\if\MODE1
\author{Saman Cyrus$^{1}$, Ali Khaki-Sedigh$^{2}$ 
	\thanks{$^{1}$ Saman Cyrus is with the Department of Electrical and Computer
		Engineering,
		University of Wisconsin--Madison and the Optimization Group at the Wisconsin
		Institute for Discovery, Madison, WI 53706, USA. 
		Email: {\tt\small cyrus2@wisc.edu}
	}%
	\thanks{$^{2}$ Ali Khaki Sedigh is with the Department of Electrical and Computer Engineering, K. N. Toosi University of Technology, Tehran, Iran.
		(e-mail: {\tt\small sedigh@kntu.ac.ir}).
		
		The authors are with the Optimization Group at the Wisconsin
		Institute for Discovery.
		Email: {\tt\small \{cyrus2,laurent.lessard\}@wisc.edu}}%
}

\markboth{IEEE Transactions on Automatic Control}%
{}
\else
\author{Saman~Cyrus\footnotemark[1] \and
	Ali~Khaki Sedigh\footnotemark[2]}
\note{}
\fi

\maketitle

\footnotetext[1]{S.~Cyrus is with the Department of Electrical and Computer
	Engineering,
	University of Wisconsin--Madison and the Optimization Group at the Wisconsin
	Institute for Discovery, Madison, WI 53706, USA. 
	Email: {\tt\small cyrus2@wisc.edu}}
\footnotetext[2]{A.~Khaki Sedigh is with the Department of Electrical and Computer Engineering, K. N. Toosi University of Technology, Tehran, Iran.
	Email: {\tt\small sedigh@kntu.ac.ir}}

\begin{abstract}
Input constrained Model predictive control (MPC) includes an optimization problem which should iteratively be solved at each time-instance. The well-known drawback of model predictive control is the computational cost of the optimization
problem. This results in restriction of the application of MPC to systems with slow dynamics, e.g., process control systems and small-scale problems. Therefore, implementing fast numerical optimization algorithms has been a point of interest. Interior-point methods are proved to be appropriate algorithms, from computational cost point-of-vie, to solve input-constrained MPC. In this paper first a modified version of
Mehrotra's predictor-corrector algorithm, a famous interior-point algorithm, is extended for quadratic programming problems and then is applied to the constrained model predictive control problems. Results show that as expected, the new algorithm is faster than Matlab solver's algorithm.
\end{abstract}

\section{Introduction}

Model predictive control (MPC) is a technique where prediction of the future behavior of a dynamical system and the resulting optimal control action are determined by iteratively solving an online optimization problem at each time-step. At each time interval, model predictive control computes a sequence of future variable adjustments by solving a finite horizon open-loop optimal control problem and using the first element of the input sequence as the control action for the
system. Then in the next iteration, all calculations are repeated. Nowadays MPC is applied to wide
range of applications from aerospace to automotive and control of chemical processes \cite{c1}.

MPC has many advantages, including its ability to be applied to multivariable systems and plants where offline
computations are not feasible. Another desired property of MPC is its capability to include explicit expression of state and input constraints.

In absence of constraints, the infinite-horizon optimization problem of MPC is the LQR problem. A well-known
drawback of MPC is its computational cost which restricts application of MPC to systems with slow
dynamics \cite{c2,c3,c6} and makes MPC inappropriate for large-scale systems \cite{c4,c11,c13}. Even for such
systems, computer technology required to perform the computations is sometimes costly and also causes software
certification concerns \cite{c12}. For such systems, often the optimal control problem in MPC is a convex
quadratic programming \cite{c2}. Different approaches have been proposed to solve the optimization problem of MPC numerically, based on linearity or nonlinearity of the system and constraints. For example, one of the simplest cases, linear-MPC, is when constraints and the dynamical system are all linear and the optimization problem's objective function is a quadratic function \cite{c10}. Another simple expression of MPC that makes it solvable is to use a hybrid prediction model, and as a result categorizing the optimization problem as a mixed integer quadratic program (MIQP) or a mixed-integer linear program (MILP) \cite{c12}.

While in general in solving the optimization problem of MPC feedback is computed online using iterative optimization algorithms, in special cases, e.g. explicit MPC \cite{c7,c19}, feedback can be computed offline. In explicit MPC, a small but dense Hessian is resulted from removing states from the constraints and the objective function. If the horizon length is shown by $N$, number of constrained variables are described by $n$, and number of inputs is shown with $p$, the computational complexity for explicit MPC is $\mathcal{O}((n+p) p^2 N^3)$, while the computational complexity of implicit MPC is shown in \cite{c2} to be $\mathcal{O}((p+m)^3+l(p+m)^2 )N)$ where $m$ is the number of states \cite{c18}. Due to memory demand, explicit MPC cannot be used for large-scale applications \cite{c7,c12,c11,c14}. Even for small problems, solving the problem online is not always slower \cite{c6}.

Increasing the speed of solving the optimization problem that we are facing with in MPC, has been a point of interest for many years. One approach is to increase the speed by parallel hardware implementation \cite{c26,c27,c28,c29}. Another approach is to increase the speed through improving numerical optimization algorithms to solve the system of equations. In the case of linear MPC, Gradient Method \cite{c25,c38,c39}, interior-point methods, and active-set methods \cite{c37,c16,c30,c35,c36} are examples of numerical algorithms that are implemented to solve MPC.

Interior-point methods (IPM) \cite{c20} have been proved as an efficient way to solve large-scale quadratic programming (QP) problems which arise in MPC \cite{c16,c18}. The reason is that interior-point algorithms can be characterized by polynomial complexity and are believed to be effective in large-scale problems, and hence MPC. Many IPMs solve Karush-Kahn-Tucker (KKT) use Newton's method \cite{c22} while inexact-IPMs do not need to solve KKT conditions \cite{c16}. 

From implementation point of view, interior-point methods and active-set methods are different. In contrast to active-set methods, interior-point methods cannot easily use warm-start strategies, although there are efforts to yield warm-start IPMs \cite{c31,c32,c33,c34}. See \cite{c15,c17,c18} for implementation of warm-start interior point methods in model predictive control context. Different interior-point methods are discussed in \cite{c22} and it has been resulted that the most practical IPM is \textit{Mehrotra's Predictor-Corrector Primal-Dual interior point method} \cite{c21,c22,c23}.

\section{Model Predictive Control}
Although the methodology of all model predictive controllers share common strategy but they may be different in system's model, objective function, and the procedure to obtain the control law. In \textit{Generalized Predictive Control} (GPC) \cite{c40}, the objective function is:
\begin{multline}
	J(N_1,N_2,N_u) \defeq \sum_{j=N_1}^{N_2} \delta (j) [\hat{y} (t+j|t)-w(t+j)]^2\\
	 + \sum_{j=1}^{N_u} \eta(j)[\Delta u(t+j-1)]^2
\end{multline}
where $N_1$ is the minimum costing horizon, $N_2$ is maximum costing horizon, $w$ is future reference trajectory, $\hat{y}(t+j|t)$ is prediction of system's output $j$ step ahead, $\delta$ and $ \eta$ are weighting sequences, $N_u$ is the control horizon and $u$ is the control sequence. Our goal is to make $w$ and $\hat{y}$ as close as possible by defining the appropriate value for $u$ sequence.

By CARIMA (Auto-Regressive Integrated Moving Average) model, system can be defined as:
\begin{equation}
	\mathscr{A}(z^{-1})y(t)=\mathscr{B}(z^{-1})z^{-d}u(t-1)+\mathscr{C}(z^{-1}) \frac{e(t)}{\Delta}
\end{equation}
where $\Delta = 1- z^{-1}$, and $\mathscr{A}$,$\mathscr{B}$, and $\mathscr{C}$ are
\begin{subequations}
	\begin{align}
	\mathscr{A}(z^{-1}) &=1+a_1 z^{-1} + a_2 z^{-2}+ \dots + a_{na}z^{-na},\\
	\mathscr{B}(z^{-1}) &=b_) + b_1 z^{-1} + b_2 z^{-2} + \dots + b_{nb} z^{-nb},\\ 
	\mathscr{C}(z^{-1}) &= 1+c_1 z^{-1} + c_2 z^{-2} + \dots + c_{nc} z^{-nc},\\
	\tilde{A}(z^{-1})   &= \Delta \mathscr{A}(z^{-1}).
	\end{align}
\end{subequations}

By Diophantine equation, we have 
\begin{equation*}
1 = E_j(z^{-1}) \tilde{A}(z^{-1}) + z^{-j}F_j(z^{-1}),
\end{equation*}
 where $j$ is the number of steps ahead. The polynomials $E_j$ and $F_j$ can be obtained by dividing $1$ by $\tilde{A}(a^{-1})$ until the remainder can be described as $z^{-j}F_j(z^{-1})$.
Having the definitions above, GPC results the following quadratic program
\begin{align}
	&J = \frac{1}{2}u^\tp G u+c^\tp u +f_0
\end{align}
where
\begin{align}
	&G=2(\Gamma^\tp \Gamma+\eta I)\\
	&c^\tp= 2(f-w)^\tp \Gamma\\
	&f_0 = (f-w)^\tp (f-w)
\end{align}
and
\begin{align}
    y &= \Gamma u+f\\
	f &= F(z^{-1})y(t)+\Gamma^{'}(z^{-1}) \Delta u(t-1)\\
	 F(z^{-1})&=\begin{bmatrix}
	F_{d+1}(z^{-1}) &\cdots&F_{d+N}(z^{-1})
	\end{bmatrix}^\tp\\
	\Gamma_j(z^{-1}) &= E_j(z^{-1}) \mathscr{B}(z^{-1})\\
	\Gamma &= \begin{bmatrix}
		g_0	&	0 &	\dots &	0\\
		g_1 & g_0 & \dots & 0\\
		\vdots & \vdots & \vdots & \vdots \\
		g_{N-1} & g_{N-2} & \dots & g_0
	\end{bmatrix}	
\end{align}
\begin{multline}
	\Gamma^{'}(z^{-1})= \\ \left[ \begin{matrix}
		\Big(G_{d+1}(z^{-1})-g_0 \Big) z\\
		\Big( G_{d+2}(z^{-1})-g_0-g_1z^{-1} \Big) z^2\\
		\vdots\\
		\Big(G_{d+N}(z^{-1})-g_0-g_1 z^{-1}-\dots - g_{N-1} z^{-(N-1)} \Big)z^N
	\end{matrix}
	\right]
\end{multline}
Here $u$, $y$, and $w$ are defined as
\begin{align*}
u &= \begin{bmatrix}
	\Delta u(t) & \cdots& \Delta u(T+N-1)	\end{bmatrix}^\tp\\
y &=\begin{bmatrix} 
	\hat{y}(t+d+1|t) &\cdots	&\hat{y}(t+d+N|t) 
\end{bmatrix}^\tp\\
w &= \begin{bmatrix}
	w(t+d+1) & \dots & w(t+d+N)
\end{bmatrix}^\tp
\end{align*}

\section{Modified Predictor Corrector Algorithm}
In \cite{c42} it is shown that Mehrotra's algorithm \cite{c21} has drawbacks. By a numerical example it is shown that there are cases where, in practice, Mehrotra's algorithm performs very small steps in order to simultaneously keep the resulting points in a specific neighborhood of the central path while also maintaining the computational complexity of the algorithm polynomial. In \cite{c42} a new algorithm with safeguard on step lengths is proposed for linear programming and it is shown that the computational complexity of the resulting algorithm is polynomial-time.

If we define the problem as
\begin{align}
	\min_x&  \quad \frac{1}{2}x^\tp G x +x^\tp c \nonumber\\
	\text{subject to}&  \quad A x \ge b 
\end{align}
where $G$ is symmetric and positive semidefinite, and matrix $A$ is $m$-by-$n$, then KKT conditions could be written as
\begin{align}
	Gx-A^\tp \lambda +c = 0& \nonumber\\
	Ax-y-b=0& \nonumber\\
	y_i \lambda_i=0& \quad i \in \{1,2,\dots,m\} \nonumber\\
	(y,\lambda) \ge 0&
\end{align}
where $y$ is the vector of slack variables and $\lambda$ is the vector of Lagrange multipliers associated with inequality constraints. Complementarity measure is defined by
\begin{equation}
	\mu = \frac{y^\tp \lambda}{m}
\end{equation}

The Mehrotra's algorithm for quadratic programming problems \cite{c23} is given in Table \ref{tab:QP Mehrotra}.

\begin{table}[!ht]
	\renewcommand{\arraystretch}{1}
	\centering
	\tabcolsep=0.11cm
	\small
	
	\begin{tabular}{l}
		\toprule
		\textbf{Mehrotra's Predictor-Corrector Algorithm}\\ \textbf{ for Quadratic Programming}\\
		\midrule
		Starting from an appropriate initial point $\displaystyle (x_0,y_0,\lambda_0)$\\
		\textbf{for} $k=0,1,2,\dots$\\[.5ex]
		\quad Set $(x,y,\lambda)=(x_k,y_k, \lambda_k)$ and solve the following\\ \quad  linear system for $(\Delta x^{\text{aff}},\Delta y^{\text{aff}}, \Delta \lambda^{\text{aff}})$\\
		\\ 
		\quad $ \left[ 
		\begin{matrix}
		G & 0 & -A^\tp\\
		A & -I & 0 &\\
		0 & \Lambda & \mathcal{Y}
		\end{matrix}
		
		\right] \left[ \begin{matrix}
			\Delta x^{\text{aff}}\\ \Delta y^{\text{aff}} \\ \Delta \lambda^{\text{aff}}
		\end{matrix}
		\right]= \left[ \begin{matrix}
			-Gx+A^\tp \lambda- c\\
			-Ax+y+b\\
			-\Lambda \mathcal{Y} e
		\end{matrix}
		\right]$
		\\
		\\
		$\quad$ Calculate $\mu = y^\tp \lambda/m$\\[.5ex]
		$\quad$ Calculate\\ 
		$\qquad\hat{\alpha}_{\text{aff}}= \text{max} \{\alpha \in (0,1] | (y,\lambda) + \alpha (\Delta y^{\text{aff}}, \Delta \lambda^{\text{aff}} ) \ge 0  \}$\\ [.5ex]
		$\quad$ Calculate \\ $\quad$ $\displaystyle \mu_{\text{aff}}=\frac{(y+\hat{\alpha}^{\text{aff}} \Delta y^{\text{aff}} )^\tp}{m} (\lambda+ \hat{\alpha}^{\text{aff}} \Delta y^{\text{aff}})^\tp(\lambda + \hat{\alpha}^{\text{aff}} \Delta \lambda^{\text{aff}})$\\[.5ex]
		$\quad$ set centering parameter to $\displaystyle \sigma = (\mu_{\text{aff}}/\mu)^3$\\[.5ex]
		$\quad$ solve the following linear system for $(\Delta x, \Delta y , \Delta \lambda)$\\[.5ex]
		\\
		$\quad \begin{bmatrix}
			G & 0 & -A^\tp\\
			A & -I & 0 &\\
			0 & \Lambda & \mathcal{Y}
		\end{bmatrix}\begin{bmatrix}
			\Delta x\\ \Delta y \\ \Delta \lambda
		\end{bmatrix}
		= \left[ \begin{matrix}
			-Gx+A^\tp \lambda- c\\
			-Ax+y+b\\
			\star
		\end{matrix}
		\right]$
		\\
		\\
		$\quad$ Choose $\tau_k \in (0,1)$ and set $\hat{\alpha}= \text{min} (\alpha^{pri}_{\tau_k},\alpha^{dual}_{\tau_k})$\\[.5ex]
		$\quad$ Set $(x_{k+1},y_{k+1},\lambda_{k+1})=(x_k,y_k,\lambda_k)+\hat{\alpha}(\Delta x, \Delta y, \Delta \lambda)$\\[.5ex]
		end \textbf{for}\\
		\bottomrule
		\end{tabular}
		\caption{\textbf{Mehrotra's Predictor-Corrector Algorithm for Quadratic Programming}}
		\label{tab:QP Mehrotra} 
		\end{table}
	
		In Table \ref{tab:QP Mehrotra}, 
		\begin{equation*}
		-\Lambda \mathcal{Y} e - \Delta \Lambda^{\text{aff}} \Delta \mathcal{Y}^{\text{aff}} e + \sigma \mu e
		\end{equation*}
		 and $\sigma \in [0,1]$. Matrices $\Lambda$ and $\mathcal{Y}$ are defined as
		\begin{align}
		\alpha_{\tau}^{\text{pri}} &= \max \{\alpha \in (0,1]: y+ \alpha \Delta y \ge (1-\tau)y \}\\
		\alpha_{\tau}^{\text{dual}} &=\max \{\alpha \in (0,1]: \lambda + \alpha \Delta \lambda \ge (1-\tau) \lambda \}
		\end{align}
		where $\tau \in (0,1)$ and
		\begin{align}
		\Lambda &= \text{diag}(\lambda_1,\dots,\lambda_m) \\
		\mathcal{Y} &= \text{diag}(y_1,y_2,\dots,y_m)\\
		e &=(1,\dots,1)^\tp.
		\end{align}
		
		The algorithm proposed by Salahi and Terlaki\cite{c42} can be extended into the quadratic programming case and the algorithm in Table \ref{tab:QP Salahi} is resulted.
		Here
		\begin{equation}
			\mathcal{N}_{\infty}^{-}(\gamma)=\{(x,\lambda,y) \in \mathcal{F}^0| \; y_i \lambda_i \ge \gamma \mu  \}
		\end{equation}
		\begin{equation}
			\gamma \in (0,1)
		\end{equation}
		\begin{equation}
			\mu = \frac{y^\tp \lambda}{m}
		\end{equation}
		\begin{equation}
			\mathcal{F}^0 =\{(x,\lambda,y) | Gx+c-A^\tp \lambda =0,\; Ax-y=b, \; (\lambda, y) \ge 0 \}
		\end{equation}
		and
		\begin{equation}
			\xi=1-\left( \frac{2 \gamma t}{1- \gamma} \right)^{\frac{1}{3}}
		\end{equation}
		\begin{equation}
			t= \max_{ i \in \mathscr{I}_{+} } \left( \frac{\Delta y_i^{\text{aff}} \Delta \lambda_i^{\text{aff}}}{y_i \lambda_i} \right)
		\end{equation}
		\begin{equation}
			\mathscr{I}= \{ 1,\dots,n \}
		\end{equation}
		\begin{equation}
			\mathscr{I}_{+}= \{i \in \mathscr{I}| \Delta y_i^{\text{aff}} \Delta \lambda_i^{\text{aff}} >0 \}
		\end{equation}
		\begin{equation}
			\mathscr{I}_{-}= \mathscr{I} \setminus \mathscr{I}_{+}
		\end{equation}

		and
		\begin{align}
			x(\alpha)=x+ \alpha \Delta x,\\ y(\alpha) = y+ \alpha \Delta y,\\ s(\alpha) = s+ \alpha \Delta s
		\end{align}

		\begin{table}[!h]
			\renewcommand{\arraystretch}{0.83}
			\raggedright
			\tabcolsep=0.1cm
			\small
			\begin{tabularx}{.4\textwidth}{l}
				\toprule
				\textbf{Revised Predictor-Corrector Algorithm for}\\ \textbf{Quadratic Programming}\\
				\midrule
				Input: parameters $\gamma \in (0,1/4)$, $\epsilon >0$ and the\\ \quad \quad \; safeguard parameter $\displaystyle \beta \in [\gamma , 1/4)$;\\ [.5ex] 
				Starting from an appropriate initial point\\ $(x_0,y_0,\lambda_0) \in \mathcal{N}_{\infty}^{-}(\gamma)$\\
				\textbf{begin}\\[.5ex]
				\quad \textbf{while} $y^\tp \lambda \ge \epsilon$ do\\
				\quad \textbf{begin}\\
				\quad \quad \textbf{Predictor Step}\\
				\quad \quad Solve the following linear system and compute \\ \quad \quad the allowable step size ($\alpha_a$)\\ \quad \quad such that $(x(\alpha_a),\lambda (\alpha_a) , y(\alpha_a)) \in \mathcal{F}$; \\ [.5ex]
				\quad \quad \quad $ \left[ 
				\begin{matrix}
				G & 0 & -A^\tp\\
				A & -I & 0 &\\
				0 & \Lambda & \mathcal{Y}
				\end{matrix}
				
				\right]  \left[ \begin{matrix}
					\Delta x^{\text{aff}}\\ \Delta y^{\text{aff}} \\ \Delta \lambda^{\text{aff}}
				\end{matrix}
				\right] = \left[ \begin{matrix}
					-Gx+A^\tp \lambda- c\\
					-Ax+y+b\\
					-\Lambda \mathcal{Y} e
				\end{matrix}
				\right]$
				\\
				\\
				\quad \quad \textbf{IF} $(1-\alpha_a)y(\alpha_a)^\tp \lambda(\xi) \le \epsilon$ \textbf{then}\\[.5ex]
				\quad \quad Let $y=y(\alpha_a), \lambda=\lambda(\alpha_a)$ and Stop.\\[.5ex]
				\quad \quad \textbf{end}\\
				\quad \textbf{end}\\
				\\
				\quad \textbf{begin}\\
				\quad \quad \textbf{Corrector Step}\\
				\quad \quad \textbf{IF} $\alpha_a > \xi$, then let $\alpha_a = \xi$\\[.5ex]
				\quad \quad \textbf{end}\\
				\quad \quad \textbf{IF} $\alpha_a \ge 0.1$, then solve the following linear system\\
				\quad \quad \; \; with $\mu=(1-\alpha_a)^3/\mu_g $ and compute the maximum\\
				\quad \quad \; \; allowable $\alpha_c$ such that $(x(\alpha_c),\lambda(\alpha_c),y(\alpha_c))\in \mathcal{N}_{\infty}^{-}$\\[.5ex]
				\quad 
				\vbox{\begin{align*}
					& \left[\begin{matrix}
					G & 0 & -A^\tp\\
					A & -I & 0 &\\
					0 & \Lambda & \mathcal{Y}
					\end{matrix}
					\right] \left[ \begin{matrix}
					\Delta x\\ \Delta y \\ \Delta \lambda
					\end{matrix}
					\right] \nonumber = \left[ \begin{matrix}
					0\\
					0\\
					\star
					\end{matrix}
					\right]
					\end{align*}}\\
				
				\quad \quad \textbf{end}\\
				\quad \quad \textbf{IF} $\alpha_a < 0.1$, then solve the following linear system\\ 
				\quad \quad \; \; with $\mu$ as defined before and compute the maximum\\ \quad \quad \; \; allowable step such that $(x(\alpha_c),\lambda(\alpha_c),y(\alpha_c))
				\in \mathcal{N}_{\infty}^{-};$\\ [.5ex]
				\quad \vbox{\begin{align*} 
						& \left[\begin{matrix}
							G & 0 & -A^\tp\\
							A & -I & 0 &\\
							0 & \Lambda & \mathcal{Y}
						\end{matrix}
						\right]  \left[ \begin{matrix}
							\Delta x\\ \Delta y \\ \Delta \lambda
						\end{matrix} \right] = \left[ \begin{matrix}
							0\\
							0\\
							\ast
						\end{matrix}
						\right] \end{align*}}\\
				\quad \quad \textbf{end}\\
				\quad \quad \textbf{IF} $\displaystyle  \alpha_c < \frac{\gamma}{\sqrt{2}n}$, then solve the following linear system\\
				\quad \quad \; \; with $ \displaystyle \mu=\frac{\beta}{1-\beta}\mu_g$ and compute the maximum\\
				\quad \quad \; \; allowable step size\\
				\quad \quad \; \; such that $(x(\alpha_c),\lambda(\alpha_c),y(\alpha_c)) \in \mathcal{N}_{\infty}^{-}$\\
				\quad \vbox{\begin{align*}
						& \left[\begin{matrix}
							G & 0 & -A^\tp\\
							A & -I & 0 &\\
							0 & \Lambda & \mathcal{Y}
						\end{matrix}
						\right] \left[ \begin{matrix}
							\Delta x\\ \Delta y \\ \Delta \lambda
						\end{matrix} \right]= \left[ \begin{matrix}
							0\\
							0\\
							\Upsilon
						\end{matrix} \right] \end{align*}}\\
				\quad \quad \textbf{end}\\
				\quad \quad Set $(x,\lambda,y)=(x(\alpha_c),\lambda(\alpha_c),y(\alpha_c))$\\
				\quad \textbf{end}\\
				\textbf{end}\\
				\bottomrule
			\end{tabularx}
		\caption{Revised Predictor-Corrector Algorithm for Quadratic Programming. $\star= -\Lambda \mathcal{Y} e - \Delta \Lambda^{\text{aff}} \Delta \mathcal{Y}^{\text{aff}} e + \sigma \mu e$ and $\ast = -\Lambda \mathcal{Y} e -\alpha_a \Delta \Lambda^{\text{aff}} \Delta \mathcal{Y}^{\text{aff}} e + \sigma \mu e$ and $\Upsilon = -\Lambda \mathcal{Y} e -\alpha_a \Delta \Lambda^{\text{aff}} \Delta \mathcal{Y}^{\text{aff}} e + \sigma \mu e$.}
		\label{tab:QP Salahi} 
		\end{table}

		\section{Simulation}
		The algorithm in Table \ref{tab:QP Salahi} is applied to the GPC problem. Different systems all including the inequality constraint ($-0.5 \le u \le 1$ ) are taken into account and comparison is made between the proposed algorithm and MATLAB's optimization toolbox version 5.0 embedded algorithm for solving quadratic programming problems using "quadprog" command which uses a reflective Newton method \cite{c44,c45}.
		
		The simulation is implemented on a computer with 3.5 GB of RAM and Intel Core2 Due Processor E7300 with clock speed 2.66 GHz. MATLAB software version 7.10.0 was used for the simulation.
		
		Results can be seen in Table \ref{Results}. Four dynamical systems with different number of poles and zeros are taken into acount. In Figure \ref{fig:1}, the output, control sequence, and  $\Delta u$ for of one of the dynamical systems is given.

	\begin{figure*}[htb]
		\begin{subfigure}{.33\textwidth}
			\centering
			\tikzsetnextfilename{MPC_1}
			\begin{tikzpicture}
			\begin{axis}[
			width=\textwidth,
			y post scale=1,
			xmin =0,xmax= 90,
			ymin =-5,ymax=10,
			ylabel={Magnitude},
			xlabel shift = -1mm,
			label style={font=\footnotesize},
			xtick={0,30,60,90},
			ytick={-5,0,5,10},
			ylabel shift= -1mm,
			legend cell align=left,
			legend style={font=\scriptsize},
			legend style={at={(1.05,0.5)},anchor=east},
			legend pos=north east,
			ticklabel style={font=\footnotesize},
			legend entries={Output Signal,Reference Signal}
			]
			\addplot [Set1-A, very thick] table [x index=0,y index=1,header=false] {Code/r.dat};
			\addplot [Set1-B,  very thick] table [x index=0,y index=1,header=false] {Code/y.dat};
			\end{axis}
			\end{tikzpicture}
			\caption{Output vs. the reference signal}
		\end{subfigure}%
		\begin{subfigure}{.33\textwidth}
			\centering
			\tikzsetnextfilename{MPC_2}
			\begin{tikzpicture}
			\begin{axis}[
			width=\textwidth,
			y post scale=1,
			xmin=0,xmax= 90,
			ymin=-2,ymax=4,
			xlabel shift = -1mm,
			ylabel shift = -4mm,
			xtick={0,30,60,90},
			ytick={-2.5,0,2.5},
			legend cell align=left,
			legend style={font=\scriptsize},
			label style={font=\footnotesize},
			legend pos=north east,
			ticklabel style={font=\footnotesize},
			legend entries={Control Signal}
			]
			\addplot [Set1-E,  very thick] table [x index=0,y index=1,header=false] {Code/ureal.dat};
			\end{axis}
			\end{tikzpicture}
			\caption{Control Signal}
		\end{subfigure}%
		\begin{subfigure}{.33\textwidth}
			\centering
			\tikzsetnextfilename{MPC_3}
			\begin{tikzpicture}
			\begin{axis}[
			width=\textwidth,
			y post scale=1,
			xmin=0,xmax= 90,
			ymin=-2,ymax=2,
			xlabel shift = -1mm,
			ylabel shift = -3mm,
			xtick={0,30,60,90},
			ytick={-1.5,0,1.5},
			legend cell align=left,
			legend style={font=\scriptsize},
			label style={font=\footnotesize},
			legend pos=north east,
			ticklabel style={font=\footnotesize},			
			legend entries={$\Delta U$}
			]
			\addplot [Set1-C,  very thick] table [x index=0,y index=1,header=false] {Code/deltaU.dat};
			\end{axis}
			\end{tikzpicture}
			\caption{Changes in the control signal}
		\end{subfigure}
		\caption{Output, control input, and $\Delta u$ for a system with $\mathscr{A} = \left[1 \; -0.8   \right]; \mathscr{B}=\left[0.4 \; 0.6 \right], \; N_u = 20 $.}
		\label{fig:1}
	\end{figure*}

%
%

		\begin{table}[!h]
			\renewcommand{\arraystretch}{1.7}
			\centering
			\tabcolsep=0.14cm
			\small
			\begin{tabular}{ c c c c c c}
				\toprule
				\textbf{$\mathscr{A}$} &\textbf{$\mathscr{B}$}  & \textbf{Alg.} & $N_u = 3$ & $N_u=10$ & $N_u = 10$ \\
				\midrule
				
				\multirow{2}{*}{\matbox{1\\ -0.8}}& \multirow{2}{*}{\matbox{ 0.4 \\ 0.6 }} & 1 & 0.5472 & 0.6684 & 0.9543\\
				 & &2 & 0.9242 &0.963545 & 1.0217\\
				\multirow{2}{*}{\matbox{1\\ -1 \\ -0.8}}& \multirow{2}{*}{\matbox{0.4 \\ 0.6}}&1 & 0.5358 & 0.6684 & 0.8579\\
				&& 2 & 0.9425 &0.9635 & 0.9785 \\
				\multirow{2}{*}{\matbox{1\\ -1 \\ -0.8}}& \multirow{2}{*}{\matbox{0.04 \\ -6}}& 1 & 0.5358 & 0.7886 & 0.8051\\
				&& 2 & 0.8961 &0.9234 & 0.8456 \\
				\multirow{2}{*}{\matbox{1\\ -1 \\ 0.675}}& \multirow{2}{*}{\matbox{0.04 \\ -6}}& 1 & 0.7489 & 0.8845 &0.8152\\
				&& 2 & 0.9124 &0.9232 & 0.9490 \\
				\bottomrule
			\end{tabular}
		\caption{Comparison between the speed of the proposed algorithm vs. MATLAB's embedded algorithm. Algorithm 1: Proposed Algorithm, Algorithm 2: MATLAB's Algorithm}
		\label{Results}
		\end{table}
		
		\section{Conclusion}
		In this paper we have proposed a revised version of Mehrotra's predictor-corrector algorithm for the optimization problem of constrained-linear model predictive control. The algorithm has been developed for 
		quadratic programming since most model predictive control algorithms solve a quadratic program to obtain the control sequence. Numerical examples show that, at least for small problems, the proposed algorithm is faster than Matlab's algorithm for solving model predictive control problems.

\bibliographystyle{abbrv}
\begin{small}
	\bibliography{references}
\end{small}

\end{document}